\documentclass[11pt]{amsart}
\usepackage[letterpaper,hmargin=1in,vmargin=1.12in]{geometry}
\usepackage{latexsym}
\usepackage{amsfonts}
\usepackage{amsthm}
\usepackage{graphicx}
\usepackage{amssymb}
\usepackage{amsmath}
\usepackage{amssymb}
\usepackage{color}

\newtheorem{remark}{Remark}

\def\classification{\@ifnextchar [{\@xfootnotenext}%
{\begingroup\let\protect\noexpand \xdef\@thefnmark{}\endgroup
\@footnotetext}}

\begin{document}


\title{Using decision problems in  public key cryptography}

\author[V.~Shpilrain]{Vladimir Shpilrain}
\address{Department of Mathematics, The City  College  of New York, New York,
NY 10031} \email{shpil@groups.sci.ccny.cuny.edu\newline \indent
http://www.sci.ccny.cuny.edu/\~{}shpil}

\author[G.~Zapata]{Gabriel Zapata}
\address{Department of Mathematics, CUNY Graduate Center, New York,
NY 10016} \email{nyzapata@verizon.net}
\thanks{Research of the first author was partially supported by
the NSF grant DMS-0405105.}

\begin{abstract}
\noindent There are several public key establishment protocols as
well as complete public key cryptosystems based on allegedly hard
problems from combinatorial (semi)group theory known by now. Most of
these problems are {\it search problems}, i.e., they are  of the
following  nature: given a property $\mathcal{P}$  and the
information that there are objects with the property $\mathcal{P}$,
find at least one particular object with the property $\mathcal{P}$.
So far, no cryptographic protocol based on a search problem in a
non-commutative (semi)group has been recognized as secure enough to
be a viable alternative to established protocols (such as RSA) based
on commutative (semi)groups, although most of these protocols are
more efficient than RSA is.

 In this paper, we suggest to use {\it decision problems} from
combinatorial group theory as the core of a public key establishment
protocol or a  public key cryptosystem. Decision problems are
problems of the following  nature: given a property $\mathcal{P}$
and an object $\mathcal{O}$, find out whether or not the object
$\mathcal{O}$ has the property $\mathcal{P}$.

 By using a popular decision problem, the {\it word problem},
we design a cryptosystem with the following features: (1) Bob
transmits to Alice an encrypted binary sequence which Alice decrypts
correctly with probability ``very close" to 1; (2) the adversary,
Eve, who is granted arbitrarily high (but fixed)  computational
speed, cannot positively identify (at least, in theory), by using a
``brute force attack", the ``1" or ``0" bits in Bob's binary
sequence. In other words: no matter what computational speed we
grant Eve at the outset, there is no guarantee that her ``brute
force attack" program will give a conclusive answer (or an answer
which is correct with overwhelming probability) about any bit in
Bob's sequence.

\end{abstract}

\maketitle

\section{Introduction}

  In search for  more efficient and/or secure alternatives to established
cryptographic protocols (such as RSA), several authors have come up
with public key establishment protocols as well as with public key
cryptosystems based on allegedly hard {\it search problems} from
combinatorial (semi)group theory, including the conjugacy  search
problem \cite{AAG, KLCHKP}, the homomorphism search problem
\cite{Grigoriev, SZ1},  the decomposition search problem
\cite{CKLHC, KLCHKP, SU}, the subgroup membership search problem
\cite{SZ2}. All these are problems of the following nature: given a
property $\mathcal{P}$  and the information that there are objects
with the property  $\mathcal{P}$, find at least one particular
object with the property  $\mathcal{P}$ from a pool $\mathcal{S}$ of
objects.

From the very  nature of these problems, one sees that security of
the corresponding cryptographic protocols relies heavily on the
assumption that the adversary has limited (typically,
subexponential) computational capabilities. Indeed, there is usually
a natural way to recursively enumerate elements of the pool
$\mathcal{S}$; the adversary can therefore just go over
$\mathcal{S}$ one element at a time until he hits one with the
property  $\mathcal{P}$ (assuming that checking the latter can be
done  efficiently). This is what is usually called the ``brute
force" attack.

In this paper, we suggest to use {\it decision problems} from
combinatorial  group theory as the core of a public key
establishment protocol or a  public key cryptosystem. Decision
problems  are problems of the following  nature: given a property
$\mathcal{P}$  and an object $\mathcal{O}$, find out whether or not
the  object  $\mathcal{O}$ has the property $\mathcal{P}$.
 Decision problems may allow us to address (to some extent) the ultimate challenge
 of public key cryptography: to
design a cryptosystem that would be secure   against (at least,
some) ``brute force" attacks by an adversary with essentially
unlimited computational capabilities.   More precisely,  our
computational model is as follows. We explain to the adversary, in
detail, how our cryptosystem works and allow him to choose, {\it up
front}, any speed of computation that he would like to have to
attack it, but after he has made his choice, he cannot change it,
i.e., he cannot accelerate his computation beyond the limit he has
chosen for himself.

 A particular decision problem that we consider here is  the {\it word problem}
which is: given a recursive presentation of a group $G$ and an
element $g \in G$, find out whether or not $g =1$ in $G$. From the
very description of the word problem we see that it consists of two
parts: ``whether" and ``not". We call them the ``yes" and ``no"
parts of the word problem, respectively. If a group is given by a
recursive presentation in terms of generators and relators, then the
``yes" part of the word problem has a recursive solution because one
can recursively enumerate all products of defining relators, their
inverses and conjugates. However, the number of factors in such a
product required to represent a word of length  $n$ which is equal
to 1 in $G$, can be very large compared to $n$; in particular, there
are groups $G$ with efficiently solvable word problem and words $w$
of length  $n$  equal to 1 in $G$, such that the number of factors
in any factorization of $w$ into a product of defining relators,
their inverses and conjugates is not bounded by any tower of
exponents in $n$, see \cite{P}. Furthermore, if in a group $G$
 the word problem is recursively unsolvable, then the length of a
 proof verifying that $w=1$ in $G$ is  not bounded by any
 recursive function of the length of $w$.

 We also note that  the ``no" part of the word problem in many groups
is recursively unsolvable, and therefore the ``brute force" attack
described above will not be effective against this part. We have to
point out though that there is no recursively presented group (or
semigroup) that would have both ``yes" and ``no" parts of the word
problem recursively unsolvable.

 Based on these general observations, we   design here a cryptographic
 protocol (see the next section) with the following features:
\begin{enumerate}

  \item  Bob
transmits to Alice an encrypted binary sequence which Alice decrypts
correctly with probability ``very close" to 1;

  \item  The adversary,
Eve, who is granted arbitrarily high (but fixed)  computational
speed, cannot positively identify (at least, in theory), by using a
``brute force attack", the ``1" or ``0" bits in Bob's binary
sequence. In other words: no matter what computational speed we
grant Eve at the outset, there is no guarantee that her ``brute
force attack" program will give a conclusive answer (or an answer
which is correct with overwhelming probability) about any bit in
Bob's sequence.
\end{enumerate}

 We note that long time ago, there was an attempt to use
 the word problem in public key cryptography \cite{MW}, but it did
 not meet with success, for several reasons. One of the reasons,
 which is relevant to the  discussion above, was pointed out quite
 recently in \cite{BMS}: the problem which is actually used in
 \cite{MW} is not the word problem, but the {\it  word  choice
 problem}: given  $g, w_1, w_2 \in G$, find out whether $g =w_1$  or   $g =w_2$ in
 $G$, provided one of the two equalities holds. In this problem,
 both parts are recursively solvable for any recursively
 presented platform group $G$ because they both are the ``yes" parts of the word
 problem, and therefore the   word  choice
 problem cannot be used for our purposes. {\color{red}  Thus, a similarity
 of our proposal to that of \cite{MW} is misleading, and ours seems
 to be the first proposal actually based on a decision problem.}

\section{The protocol}
\label{protocol}

  Here is a sketch of our cryptographic protocol;
   details are given in the following sections.

\medskip \emph{\textbf{Protocol:}}

\begin{enumerate}

    \item A  pool of group presentations with efficiently
 solvable word problem
is considered public (e.g. is part of Alice's software).

\item Alice chooses randomly a particular presentation $\Gamma$ from the pool,
diffuses it by isomorphism-preserving  transformations to obtain a diffused presentation
$\Gamma'$, discards some of the relators and publishes the abridged diffused presentation ${\hat
\Gamma}$.

\item Bob transmits  his private binary sequence to Alice by transmitting an
element equal to 1 in ${\hat \Gamma}$ (and therefore also in
$\Gamma'$) in place of ``1" and an element not equal to 1 in
$\Gamma'$ in place of ``0".

\item Alice recovers Bob's binary sequence by first converting
elements of $\Gamma'$ to the corresponding (under the isomorphism
that she knows) elements of $\Gamma$, and then solving the word
problem in $\Gamma$.

\end{enumerate}

 Most parts of this protocol are rather nontrivial and open several
interesting research avenues. We discuss parts (1), (2), (3) in our
Sections \ref{pool}, \ref{Tietze}, \ref{equal?}, respectively.

{\it A priori} it looks like the most nontrivial part is finding an
element which is not equal to 1 in
 $\Gamma'$ since Bob does not even know the whole presentation
$\Gamma'$. We solve this problem by ``going with the flow",
 so to speak. More specifically, we just let Bob select a random
 (well, almost random)  word of sufficiently big length and show that, with overwhelming
 probability, such an element is not equal to 1 in
 $\Gamma'$. We discuss this in more detail
in Section \ref{equal?}.

 We emphasize once again what is, in our
opinion, the main advantage (at least,  theoretical) of our protocol
over the existing ones. The point is to deprive the adversary (Eve)
from attacking the protocol by doing an ``exhaustive search", which
is the most obvious (although, perhaps, often ``computationally
infeasible") way to attack all existing public key protocol.

 The way we plan to achieve our goal is relevant to part (3) of
the above protocol, more specifically, to solving the word problem
in ${\hat \Gamma}$. If Bob transmits an element $g$ equal to 1 in
${\hat \Gamma}$, Eve may be able to detect this by going over all
products of all conjugates of relators from ${\hat \Gamma}$ and
their inverses. This set is recursive, but as we have pointed out in
the Introduction, there are groups $G$ with efficiently solvable
word problem and words $w$ of length  $n$  equal to 1 in $G$, such
that the length of a proof verifying that $w=1$ in $G$ is  not
bounded by any tower of exponents in $n$, see \cite{P}.

Furthermore, if Bob transmits an element $g$ not equal to 1 in
${\hat \Gamma}$, then detecting this is even more difficult for Eve.
In fact, it is impossible in general; Eve's only hope here is that
she will be lucky to find a factor group of ${\hat \Gamma}$ where
the word problem is solvable, and that $g \ne 1$ in that factor
group. This is what we call a {\it  quotient attack}, see our
Section \ref{quotient}.

\begin{remark} It may look like an encryption
 protocol  with the features outlined in the Introduction cannot exist; in particular,
the following attack by a computationally superior  adversary, Eve,
may seem viable:

\begin{quote}
 Eve can perform key generations   over and over again,
each time with fresh randomness, until the public key to be attacked
is obtained -- this will happen eventually with overwhelming
probability. Already the correctness (no matter if perfect or only
with overwhelming probability) of the scheme guarantees  that the
corresponding secret key (as obtained by Eve while  performing key
generation) allows to decrypt illegitimately.
\end{quote}

 This would be indeed viable if the correctness of the legitimate decryption by Alice was
perfect. However, in our situation this kind of attack will not work
for a general ${\hat \Gamma}$. Suppose Eve is building up two lists,
corresponding to two possible encryptions ``~$0 \to w \ne 1$ in
${\hat \Gamma}$" or ``$1 \to w=1$ in ${\hat \Gamma}$" by Bob. Our
first observation is that the list that corresponds to
 ``~$0 \to w \ne 1$" is useless to Eve because it is simply going to
 contain {\rm all} words in the alphabet $X = \{x_1, \dots, x_n, x_1^{-1},
\dots, x_n^{-1}\}$ since Bob is choosing such $w$ simply as a random
  word. Therefore, Eve may just as well forget about this list and concentrate on the other
one, that corresponds to ``$1 \to  w=1$".

Now the situation boils down to the following: if a  word  $w$
 transmitted by Bob appears on the list, then it is equal to 1 in
 $G$. If not, then not. The only problem is: how can Eve possibly
 conclude that $w$ does {\rm  not} appear on the list if the list is
 infinite? Our opponent could say here that Eve can stop at some
 point and conclude that $w \ne 1$ with overwhelming
probability, just like Alice does. The point  however is  that this
probability may not at all be as ``overwhelming" as the probability
of the  correct decryption by Alice. Compare:

\begin{enumerate}

    \item For Alice to decrypt correctly ``with overwhelming
probability", the probability $P_1(N)$
    for a random word  $w$ of length  $N$ not to be equal to 1
    should converge to 1 (reasonably fast) as $N$ goes to infinity.

\item For Eve  to decrypt correctly ``with overwhelming
probability", the probability $P_2(N, f(N))$ for a random word  $w$
of length  $N$, which is equal to 1, to have a {\rm  proof} of
length $\le f(N)$ verifying that $w=1$, should converge to 1
(reasonably fast) as $N$ goes to infinity. Here $f(N)$ represents
Eve's computational capabilities; this function can be arbitrary,
but fixed.

\end{enumerate}

 We  see that the functions $P_1(N)$ and  $P_2(N)$ are of very
 different nature, and any correlation between them is unlikely.
We note that the function $P_1(N)$ is generally well understood, and
in particular, it is known that in any infinite group $G$, $P_1(N)$
indeed converges to 1 as $N$ goes to infinity; see our Section
\ref{equal?} for more details.

 On the other hand, the functions $P_2(N, f(N))$ are more complex; they are currently
subject of a very active research, and  in particular, it appears
likely that for any $f(N)$, there are groups in which $P_2(N, f(N))$
does not converge to 1 at all. Of course, $P_2(N, f(N))$ may depend
on a particular algorithm used by Bob to produce words equal to 1,
but we leave this  discussion to  another paper.

We also note, in passing, that  if in a group $G$
 the word problem is recursively unsolvable, then the length of a
 proof verifying that $w=1$ in $G$ is  not bounded by any
 recursive function of the length of $w$.
\end{remark}

 To conclude this section, we guide the reader to
 other sections of this paper where the questions of efficiency (for legitimate parties) are
 addressed. All  steps of our protocol
 are shown to be quite efficient with the suggested
 parameters (the latter are summarized in   Section \ref{parameters}).
 Step (2) of the protocol (Alice's algorithm for obtaining her public and private keys) is
 discussed in Section \ref{Tietze}. Step (3) (encryption by Bob) is
 discussed in Section \ref{equal?}. It turns out that encryption (of
one bit) takes quadratic time in the length of a transmitted word;
the latter is approximately 150   on average, according to our
computer experiments. Step (4) (decryption by Alice)
 is discussed at the end of Section \ref{Tietze}. It is straightforward to see that
 the time Alice needs to decrypt each transmitted word $w$ is bounded
 by $C \cdot |w|$, where $|w|$ is the length of $w$ and $C$ is a
 constant which basically depends on Alice's private isomorphism between $\Gamma$  and
 $\Gamma'$.

 The fact that  Alice (the receiver) and  the adversary are separated in
 power is essentially due to Alice's knowledge of her private isomorphism
 between $\Gamma$  and  $\Gamma'$ (note that Bob does {\it not} have
 to know this isomorphism for encryption!).

 We have to admit here one disadvantage of our protocol compared to
 most well-established public key protocols: we have encryption with a rather big
 ``expansion factor". Computer experiments show that, with suggested
 parameters, one bit in Bob's message gets encrypted into a word
 of  length approximately 150   on average. This is the price we have to
 pay for granting the adversary too much computational power.

Finally, we  touch upon semantic security in the   end of Section
\ref{equal?}.

\section{Pool of group presentations}
\label{pool}

There are many classes of finitely presented groups with solvable
word problem known by now, e.g.
 one-relator groups, hyperbolic groups, nilpotent groups, metabelian groups. Note however that
 Alice should be able to randomly select a presentation from the pool {\it efficiently}, which imposes
 some restrictions on classes of presentations that can be used in this context.
 The class of finitely presented groups  that we suggest to include in our pool
 is {\it small cancellation groups}.

Small cancellation groups have relators satisfying a simple (and
efficiently verifiable) ``metric condition" (we follow the
exposition in \cite{L-S}). More specifically, let $F(X)$ be the free
group with a basis $X = \{\, x_i \,|\, i\in I \,\}$, where $I$ is an
indexing set. Let $\epsilon_k\in \{\pm 1\}$, where $1\leq k\leq n$.
A word
$w(x_1,\ldots,x_n)=x_{i_1}^{\epsilon_1}x_{i_2}^{\epsilon_2}\cdots
x_{i_n}^{\epsilon_n}$ in $F(X)$, with all   $x_{i_k}$ not
necessarily distinct, is a \emph{reduced $X$-word} if
$x_{i_k}^{\epsilon_k}\neq x_{i_{k+1}}^{-\,\epsilon_{k+1}}$, for
$1\leq k \leq n-1$. In addition, the word $w(x_1,\ldots,x_n)$ is
\emph{cyclically reduced} if it is a reduced $X$-word and
$x_{i_1}^{\epsilon_1} \neq x_{i_{n}}^{-\,\epsilon_n}$. A set $R$
containing cyclically reduced words from $F(X)$ is
\emph{symmetrized} if it is closed under cyclic permutations and
taking inverses.

 Let $G$ be a group  with presentation $\langle X;R\rangle$. A
non-empty word $u \in F(X)$ is called a \emph{piece} if there are
two distinct relators  $r_1, r_2 \in R$ of $G$ such that $r_1 = u
v_1$ and $r_2 = u v_2$ for some $v_1, v_2 \in F(X)$, with no
cancellation between $u$ and  $v_1$  or between $u$ and  $v_2$. The
group $G$  belongs to the class
  $C(p)$   if no element of $R$ is a product of fewer than $p$ pieces. Also, the group
  $G$ belongs to the class $C'(\lambda)$
if  for every $r\in R$ such that $r = uv$ and $u$ is a piece,  one
has  $|u| < \lambda |r|$.

 In particular, if $G$ belongs to the class    $C'(\frac{1}{6})$,
then Dehn's algorithm solves the word problem  for $G$ efficiently.
This algorithm is very simple: in a given word $w$, look for a
``large" piece of a relator from $R$ (that means, a piece whose
length is more than a half of the length of the whole relator). If
no such piece exists, then $w \ne 1$ in $G$. If such a piece, call
it $u$, does exist, then $r= uv$ for some  $r \in R$, where the
length of $v$ is smaller than that of $u$. Then replace $u$ by
$v^{-1}$ in $w$. The length of the resulting word is smaller than
that of $w$; therefore, the algorithm will terminate in a finite
number of steps. It has quadratic time complexity with respect to
the length of $w$.

 We also note that a generic finitely presented group is a small cancellation group (see
\cite{AO});
 therefore, to randomly select a small cancellation group, Alice can just take a few random
 words and check whether the  corresponding symmetrized set satisfies the condition for
 $C'(\frac{1}{6})$. If not, then repeat.

 To conclude this section, we give a more specific recipe for Alice
 to produce a presentation $\Gamma$  for the protocol in Section
 \ref{protocol}.

\begin{enumerate}

\item  Alice fixes a number  $k, ~10 \le k \le 20$, of generators
in her presentation $\Gamma$.  Her  $\Gamma$ will therefore have
generators  $x_1, \ldots, x_k$.

\item  Alice selects  $m$ random words  $r_1, \ldots, r_m$ in the
generators $x_1, \ldots, x_k$. Here  $10 \le  m \le 30$ and the
lengths $l_i$ of $r_i$ are random integers from the interval $L_1
\le  l_i \le L_2$. Particular values that we suggest are: $L_1 = 12,
~L_2 =20$.

\item After Alice obtains the abridged  presentation ${\hat
\Gamma}$, she adds a relation
$$x_i' = \prod_{j=1}^{M} [x_i', w_j]$$
 to it, where $x_i'$ is a (randomly chosen) generator from ${\hat \Gamma}$,   $~w_j$ are
random elements of length 1 or 2 in the generators $x_1', x_2',
\ldots, $ and  $M=10$. (Our commutator notation is:
$[a,b]=a^{-1}b^{-1}ab$.)   This relation  is needed to foil {\it
quotient attacks}, see Section \ref{quotient}. Then Alice finds the
preimage of this relation under the isomorphism between $\Gamma$ and
${\hat \Gamma}$ and  adds this preimage to the defining relators of
$\Gamma$. Thus, $\Gamma$ finally has $k$ generators and $m+1$
defining relators.

\item Finally, Alice checks whether her private presentation
$\Gamma$   satisfies the small cancellation condition
$C'(\frac{1}{6})$ (it  will with overwhelming probability, see
\cite{AO}). If not, then she has to start over.

\end{enumerate}

\section{Tietze transformations: elementary isomorphisms}
\label{Tietze}

 In this section, we explain how Alice can implement step (2) of the protocol given in the
 Introduction. First we introduce Tietze transformations; these are ``elementary isomorphisms":
 any isomorphism between finitely presented groups is a composition of Tietze transformations.
What is important to us is that every Tietze transformation is easily invertible, and therefore
Alice can compute the inverse isomorphism that takes $\Gamma'$ to $\Gamma$.

 Tietze introduced isomorphism-preserving elementary transformations that can be applied to
groups presented by generators and relators. They are of the
following  types.
\begin{description}
\item[(T1)]
\emph{Introducing a~new generator}: Replace $\langle x_1,x_2,\ldots
\mid r_1, r_2,\dots\rangle$ by\\
 $\langle y, x_1,x_2,\ldots \mid
ys^{-1}, r_1, r_2,\dots\rangle$, where $s=s(x_1,x_2,\dots )$ is an
arbitrary element in the generators $x_1,x_2,\dots$.
\item[(T2)]
\emph{Canceling a~generator} (this is the converse of (T1)): If we have a~presentation of the form
$\langle y, x_1,x_2,\ldots \mid q, r_1, r_2,\dots\rangle$, where $q$ is of the form $ys^{-1}$, and
$s, r_1, r_2,\dots $ are in the group generated by $x_1,x_2,\dots$, replace this presentation by
$\langle x_1,x_2,\ldots \mid r_1, r_2,\dots\rangle$.
\item[(T3)]
\sloppy \emph{Applying an    automorphism}: Apply an automorphism of
the free group generated by $x_1,x_2,\dots $ to all the relators
$r_1, r_2,\dots$.
\item[(T4)]
\emph{Changing defining relators}: Replace the set $r_1, r_2,\dots$
of defining relators by another set $r_1', r_2',\dots$ with the same
normal closure. That means, each of  $r_1', r_2',\dots$ should
belong to the normal subgroup generated by $r_1, r_2,\dots$, and
vice versa.
\end{description}

Tietze has proved (see e.g. \cite{L-S}) that two groups $\langle
x_1,x_2,\ldots \mid r_1, r_2,\dots\rangle$ and $\langle
x_1,x_2,\ldots \mid s_1, s_2,\dots\rangle$ are isomorphic if and
only if one can get from one of the presentations to the other by a
sequence of transformations \textup{(T1)--(T4)}.

 For each Tietze transformation of the types \textup{(T1)--(T3)}, it is easy to obtain an explicit
 isomorphism (as a mapping on
 generators) and its inverse. For a Tietze transformation of the
 type \textup{(T4)}, the isomorphism is just the identity map.
We would like here to make Tietze transformations of the
 type \textup{(T4)} recursive, because {\it a priori} it is not clear how
 Alice can actually apply these transformations. Thus, Alice is
 going to use the following recursive version of \textup{(T4)}:

\medskip

\noindent {\bf (T4$'$)} In the set $r_1, r_2,\dots$, replace some
$r_i$ by one of the:  $r_i^{-1}$,  $r_i r_j$, $r_i r_j^{-1}$, $r_j
r_i$, $r_j r_i^{-1}$, $x_k^{-1} r_i x_k$, $x_k r_i x_k^{-1}$, where
$j \ne i$, and $k$ is arbitrary.
\medskip

 We  suggest  that in part (2) of the protocol in Section
 \ref{protocol}, Alice should first apply several transformations
 of the type \textup{(T4$'$)} to ``mix" the presentation $\Gamma$.
(This does not add complexity to the final isomorphism since  for a
Tietze transformation of the
 type \textup{(T4)}, the isomorphism is just the identity map, as we have noted
 above.)
 In particular, if $\Gamma$ was a small cancellation presentation (see Section
 \ref{pool}) to begin with, then after applying several
 transformations
 \textup{(T4$'$)} it will, most likely, no longer be. As a result,
 Eve's chances to augment the public presentation ${\hat \Gamma}$ to
a small cancellation presentation (see Section
 \ref{isomorphism})  are getting slimmer.

One more trick that Alice can use for better diffusion of her
presentation  is making a free product of her group with the trivial
group given by a non-standard presentation. That means, she can add
new generators $z_1, \ldots, z_q$ and new relators  $s_1(z_1,
\ldots, z_q), \ldots, s_t(z_1, \ldots, z_q)$, such that the
presentation $\langle z_1, \ldots, z_q \mid s_1,  \ldots, s_t
\rangle$ defines the trivial group. After that, she has to apply
several (T3)s and (T4$'$)s to mix the new generators with the old
ones. We note that there are many non-trivial presentations of the
trivial group to choose from; for example, in \cite{MMS}, there are
given several infinite series of such presentations in the special
case where $t=q$ (so-called {\it balanced presentations}). Without
this restriction, there are even more choices; in particular, Alice
can just add arbitrary relators to a balanced presentation of the
trivial group, thus adding to the confusion of the adversary.

After Alice has mixed $\Gamma$ by using these tricks, we suggest
that she should aim for  breaking down some of the defining relators
into ``small pieces". More formally, she can replace a given
presentation by an isomorphic presentation  where
 most defining relators have length at most 4. (Intuitively, diffusion
of elements should be easier to achieve in a group with  shorter
defining relators). This is easily achieved by applying
transformations \textup{(T1)}  (see below) which can be ``seasoned"
by a few elementary automorphisms (type \textup{(T3)}) of the form
$x_i \to x_i x_j^{\pm 1}$  or  $x_i \to x_j^{\pm 1}x_i $, for better
diffusion.

 The procedure of breaking down defining relators is quite simple. Let $\Gamma$ be a presentation
$\langle x_1,\dots , x_k; r_1,...,r_m\rangle$. We are going to
obtain a different, isomorphic, presentation by using Tietze
transformations of types \textup{(T1)}.
 Specifically, let, say,  $r_1=x_i x_j u$, $1 \le i,j \le k$.
We introduce a new generator $x_{k+1}$  and a new relator
$r_{m+1}=x_{k+1}^{-1}x_i x_j$. The presentation \\
$\langle x_1,\dots
,x_k, x_{k+1}; r_1,\dots ,r_m, r_{m+1}\rangle$ is obviously
isomorphic to $\Gamma$. Now if we replace $r_1$ with $r_1'=x_{k+1}
u$, then the presentation $\langle x_1,\dots ,x_k, x_{k+1};
r_1',\dots ,r_m, r_{m+1}\rangle$ will again be isomorphic to
$\Gamma$, but now the length of one of the defining relators ($r_1$)
has decreased by 1. Continuing in this manner, Alice can eventually
obtain a presentation where many relators have length at most 3, at
the expense of introducing more generators. In fact, relators of
length 4 are  also good for the purpose of diffusing a given word,
so we are not going to ``cut" the relators into too small pieces
(i.e., we do not want pieces of length 1 or 2), but rather settle
with relators of length 3 or 4. Most of the longer relators can be
discarded from the presentation $\Gamma'$ to obtain the abridged
  presentation ${\hat \Gamma}$.

 We conclude this section with a simple  example, just to illustrate
 how Tietze transformations can be used to cut relators into
 pieces. In this example, we start with a presentation having two relators of
length 5 in 3 generators, and end up with a presentation having 4
relators of length 3 or 4 in 5 generators. The symbol $\cong$ below
means ``is isomorphic to".

\medskip

\noindent {\bf Example.}  $\langle x_1, x_2, x_3 ~\mid ~x_1^2x_2^3,
~x_1x_2^2x_1^{-1}x_3 \rangle ~\cong ~\langle x_1, x_2, x_3, x_4
~\mid ~x_4=x_1^2, ~x_4x_2^3, ~x_1x_2^2x_1^{-1}x_3 \rangle\\
\cong \langle x_1, x_2, x_3, x_4, x_5 ~\mid  ~x_5=x_1x_2^2,
~x_4=x_1^2, ~x_4x_2^3, ~x_5x_1^{-1}x_3 \rangle ~\cong \\
\langle x_1, x_2, x_3, x_4, x_5 ~\mid  ~x_5=x_2^2, ~x_4=x_1^2,
~x_4x_2^3, ~x_1x_5x_1^{-1}x_3 \rangle$.
\medskip

The last isomorphism illustrates applying a transformation of type
(T3), namely, the automorphism  $x_5 \to x_1x_5, ~x_i \to x_i, ~i
\ne 5$.

\section{Generating random elements in finitely presented groups}
\label{equal?}

 In this section, we explain how to  implement the crucial step (3) of the protocol given in
Section \ref{protocol}.

 When Bob wants to transmit an
element equal to 1 in ${\hat \Gamma}$, he should construct a word,
looking ``as random as possible" (for semantic security), in the
relators ${\hat r_1},\dots ,{\hat r_l}$ and their conjugates. When
he wants to transmit an element not equal to 1 in ${\hat \Gamma}$,
he just
 selects a random word of sufficiently big length; it turns out that, with overwhelming
 probability, such an element is not equal to 1 in
 $\Gamma'$ (we explain it in the end of this section).

We start with a description of Bob's possible diffusion strategy for
producing elements equal to 1 in ${\hat \Gamma}$. (It is rather
straightforward to produce a random word of a given length in
generators $x_1,\dots ,x_k$, so we are not going to discuss it
here.) When Bob transmits a word $w$ equal to 1 in ${\hat \Gamma}$,
he wants to diffuse it so that large pieces of defining relators
would not be visible in $w$. In some specific groups (e.g. in braid
groups) a diffusion is provided by a ``normal form", which is a
collection of symbols that uniquely corresponds to a given element
of the group. The existence of such normal forms is usually due to
some special algebraic or geometric properties of a given group.

However, since Bob does not know any meaningful properties of the
group defined by the presentation ${\hat \Gamma}$ which is given to
him, he cannot employ normal forms in the usual sense. The only
useful property that the presentation ${\hat \Gamma}$ has is that
most of its defining relators  have length 3 or  4, see Section
\ref{pool}. We are going to take advantage of this property as
follows. We suggest the following procedure which is probably best
described by the word ``shuffling".


\begin{enumerate}

  \item Make a   product of the form $u=s_1 \cdots s_p$, where
  each $s_i$ is randomly chosen among defining relators ${\hat r_1}, {\hat r_2},
  \dots$, of length 3 or  4,
  their inverses, and their
  conjugates by one- or two-letter words in $x_1', x_2', \dots $. The
  number $p$  of factors should be sufficiently big, at least 10 times
  the number of defining relators in ${\hat \Gamma}$.

\item Insert approximately $\frac{2p}{k}$ expressions of the form
$x_j' (x_j')^{-1}$  or  $(x_j')^{-1}x_j'$ in random places of the
word $u$ (here $k$ is the number of generators  $x_i'$ of ${\hat
\Gamma}$), for random values of $j$.

\item Going left to right in the word $u$, look for two-letter
subwords that are parts of  defining relators ${\hat r_i}$ of length
3 or 4. When you spot such a subword, replace it by the inverse of
the augmenting part of the same defining relator and continue. For
example, suppose there is a relator ${\hat r_i}=x_1'x_2'x_3'x_4'$,
and suppose you spot the subword $x_1'x_2'$ in $u$. Then replace it
by $(x_4')^{-1}(x_3')^{-1}$ (obviously,
$x_1'x_2'=(x_4')^{-1}(x_3')^{-1}$ in your group). If you spot the
subword $x_2'x_3'$ in $w$, replace it by $(x_1')^{-1}(x_4')^{-1}$.
If there is more than one choice for replacement, choose randomly
between them.

\item Cancel remaining subwords (if any) of the form
$x_j' (x_j')^{-1}$  or  $(x_j')^{-1}x_j'$.


\end{enumerate}

Steps (2)--(4)  should be repeated approximately $p$ times for good
 mixing.

  Finally, after Bob has obtained a word $u$ this way, he sets
  $w=[x_i', u]$ and applies steps (2)--(4) to $w$ approximately $\frac{|w|}{2}$ times,
  where $|w|$ is the length of $w$. This final
  step is needed to make this $w$ (which is equal to 1 in ${\hat
  \Gamma}$, and therefore also in $\Gamma'$) indistinguishable from $w \ne 1$,
which is constructed
  in the same form $w=[x_i', u]$, see below. Here $x_i'$ is the same as in the relator
$x_i' = \prod_{j=1}^{M} [x_i', w_j]$ published by Alice, see Section
\ref{pool}. Having  $w \ne 1$ in this form is needed,   in turn,
to foil ``quotient attacks", see the end of Section \ref{quotient}.

 When Bob wants to transmit an element not equal to 1 in
$\Gamma'$, he should first choose a random word $u$ from the
commutator subgroup of the free group generated by $x_1', x_2',
\ldots $. To select a  random word from the commutator subgroup is
easy; Bob can select an arbitrary random word $v$ first, and  then
adjust the exponents on the generators in $v$ so that the exponent
sum on every  generator in $v$ is 0. The length of $u$ should be in
the same range as the lengths of the words $u$ equal to 1 in ${\hat
  \Gamma}$  constructed by Bob before. Then Bob lets $w=[x_i', u]$,
where $x_i'$ is the same as in the relator  $x_i' = \prod_{j=1}^{M}
[x_i', w_j]$ published by Alice, see Section \ref{pool}. Finally, to
hide $u$, he applies ``shuffling" to $w$ (steps (2)--(4) above)
approximately $\frac{|w|}{2}$ times, where $|w|$ is the length of
$w$.

\medskip

 Now we explain why a random word of sufficiently big
 length is not equal to 1 in $\Gamma$ with overwhelming
 probability, provided $\Gamma$  is a presentation
described in the end of Section \ref{pool}.

 Like any other group, the group $G$ given by the presentation $\Gamma$  is
 a factor group $G=F/R$ of the ambient free group $F$ generated by  $x_1,
 x_2, \dots $. Therefore, to estimate the probability that a random word in
$x_1, x_2, \dots $ would not belong to $R$ (and therefore, would not
 be equal to 1 in $G$), one should estimate  the {\it asymptotic
density} (see e.g. \cite{KMSS}) of the complement to $R$ in the free
group $F$. It makes notation simpler if one deals instead with the
asymptotic density of $R$ itself, which is

$$\rho_{_{\tiny F}}(R) = \limsup_{n\to\infty}\frac{\#\{u\in R: l(u)\le
  n\}}{\#\{u\in F: l(u)\le n\}}.$$

  Here $l(u)$ denotes the usual lexicographic length of $u$ as a word in
$x_1, x_2, \dots $.             Thus, the asymptotic density
depends, in general, on a free generating set of $F$, but we will
not go into these details here because all facts that we are going
to need are independent of the choice of basis. One  principal fact
that we can use here is due to Woess  \cite{Woess}:  if the group
$G=F/R$ is infinite, then $\rho_{_{\tiny F}}(R) = 0$. Since the
group $G$ given by the presentation $\Gamma$  is infinite (see our
Section \ref{pool}), this already tells us that the probability for
a random word of length $n$ in $x_1, x_2, \dots $  not to be equal
to 1 in $G$ is approaching 1 when $n\to\infty$. However, if we want
words transmitted by Bob to be of reasonable length (on the order of
100--200, say), we have to address the question of {\it how fast}
the ratio in the definition of the asymptotic density converges to 0
if $R$ is the normal closure of the relators described in the end of
Section \ref{pool}. It turns out that for {\it non-amenable} groups
the convergence is exponentially fast; this is also due to Woess
\cite{Woess}. We are not going to explain here what amenable groups
are; it is sufficient for us to know that small cancellation groups
are not amenable (because they have free subgroups, see e.g.
\cite{L-S}). Thus, small cancellation groups are just fine for our
purposes here: the probability for a random word of length $n$ in
$x_1, x_2, \dots $  not to be equal to 1 in $G$ is approaching 1
exponentially fast when $n\to\infty$.

Finally, we  touch upon semantic security (see \cite{GM}) of the
words transmitted by Bob. We do not give any rigorous probabilistic
estimates since this would require at least defining  a probability
measure on an infinite group, which is a very nontrivial problem by
itself (cf. \cite{BMSprobab}). Instead, we offer here an informal
argument  which we hope to be convincing. A nice thing about Bob's
encryption procedure is that when he selects a word $u \ne 1$, he
simply selects a random word. Thus, $u \ne 1$ is indistinguishable
from a random word just because it is random! Then, the element
$w=[x_i', u]$, which is   transmitted by Bob, looks like it is no
longer random because it is of a special form. However:

\begin{enumerate}

\item What is actually transmitted by Bob is a {\it word in the
alphabet $x_1', x_2', \ldots$} representing the element $w=[x_i',
u]$ of the group defined by ${\hat \Gamma}$. This word is {\it not}
of the form $[x_i', u]$ because Bob has applied a ``shuffling" to
$w$.

\item Given the specifics of our protocol, what really matters is that
transmitted  words equal to 1 in ${\hat \Gamma}$ are
indistinguishable from transmitted  words not equal to 1. This is
why we require Bob's elements representing 1 in ${\hat \Gamma}$ to
be of the form $[x_i', u]$ as well.

\end{enumerate}

 Thus, the question about semantic security of Bob's transmissions
 boils down to the following question of independent interest:
  is a word $u$ representing 1 in ${\hat
 \Gamma}$ indistinguishable from a random word (of  the
 same length)? As we have admitted above, we do not have a rigorous
 proof that it is, but computer experiments show that when most of
 the relators in ${\hat \Gamma}$ have length at most 4, then the  words $u$ representing 1 in ${\hat
 \Gamma}$, obtained as described  earlier in this section, pass at
 least the equal frequency test for 1-, 2-, and  3-letter subwords,
 thus making it appear likely that the answer to the question above
 is affirmative for such ${\hat \Gamma}$.

\section{Suggested parameters}
\label{parameters}

 In this section, we summarize all suggested parameters of our
 protocol for the reader's convenience, although most of these
 parameters were already discussed in previous sections.

\begin{enumerate}

\item The number of generators $x_i$ in Alice's private presentation
$\Gamma$ is $k$, a random integer from the interval  $10 \le k \le
20$.

\item Relators $r_1, \ldots, r_m$ in Alice's private presentation
$\Gamma$ are random words  in the generators $x_1, \ldots, x_k$.
Here $m$ is  a random integer from the interval  $10 \le m \le 30$,
and the lengths $l_i$ of $r_i$ are random integers from the interval
$12 \le  l_i \le 20$. There is one other, special, relator in
$\Gamma$,
 which is obtained as described in the end of Section \ref{pool}.

\item Alice's private isomorphism (between the presentations
$\Gamma$ and  $\Gamma'$) is a product, in random order, of $s_1$
elementary transformations of type (T1) and $s_2$  elementary
transformations of type (T3) (see Section \ref{Tietze}). Each
relator introduced by a transformation of type (T1) is  a random
word whose length is  a random integer from the interval  $[12,
20]$. Neither of the parameters $s_1, s_2$ is specified, but their
sum should be at least 50; more specifically,  as soon as $s_1+s_2$
becomes equal to 50, only transformations of type (T1) are applied,
targeted at making about 30\% of the relators to have length at most
4, as described in Section \ref{Tietze}. After that, Alice should
discard about 70\% of the relators, taking care that among the
remaining relators, at least 50\% have length at most 4.

\item Bob encrypts his secret bits by words in the given alphabet,
as described in Section \ref{equal?}. Here we specify the length of
those words. Recall that Bob starts building a   word  $w=1$ in
$\Gamma'$ as a product of $p$ words randomly chosen among published
defining relators, their inverses, and their
  conjugates by one- or two-letter words in the published
  generators. We specify $p$ as a random integer from the interval  $5
  \le p \le  12$, thus making the whole $w$ a word of length
  approximately 150 on average. Computer experiments show that
  subsequent ``shuffling", as described in Section \ref{equal?},
  only slightly increases the length of $w$.

 Finally, we recall that Bob selects a word  $w \ne 1$ in $\Gamma'$
in the form $[x_i', u]$, where $u$  is a  random word from the
commutator subgroup of the free group generated by the published
  generators $x_1', x_2', \ldots $. We therefore specify $u$ as a
  random word of length $l$, where $l$  is a  random integer from the interval
$65 \le l \le 85$. Then the length of $w=[x_i', u]$, which is
$2l+2$, is going to be approximately 150 on average, just as in the
case of $w=1$ considered above.

\end{enumerate}

\section{Isomorphism attack}
\label{isomorphism}

In this section, we discuss a (theoretically) possible ``brute
force" attack on the protocol from Section \ref{protocol}.

Knowing the pool of group presentations from which Alice selects her
private presentation $\Gamma$, Eve can try to augment the public
presentation ${\hat \Gamma}$ to a  presentation that would be
isomorphic to one from the pool. Theoretically, this is possible
because the pool is recursive and because the set of finite
presentations isomorphic to a given one is recursive, too. However,
this procedure requires enormous resources. Let us take a closer
look at it.

Eve can add to ${\hat \Gamma}$  one element at a time and check
whether the  resulting presentation, call it ${\hat \Gamma}_+$, is
isomorphic to one of the presentations from Alice's pool. The latter
is done the following way. Suppose Eve wants to check whether ${\hat
\Gamma}_+$ is isomorphic to some $\Gamma_i$. She goes over mappings
from $\Gamma_i$ to ${\hat \Gamma}_+$, one at a time, defined on the
generators of $\Gamma_i$. At the same time, she also goes over
mappings from ${\hat \Gamma}_+$ to  $\Gamma_i$ defined on the
generators of ${\hat \Gamma}_+$. She composes various pairs of these
mappings and checks: (1) whether she gets the identical mapping on
$\Gamma_i$,  and   (2) whether both mappings in such a pair are
homomorphisms, i.e., whether  they send relators of either
presentation to elements equal to 1 in the other presentation.
Having the word problem in $\Gamma_i$ solvable makes the former
checking more efficient, but it is, in fact, not necessary because
what matters here is the ``yes" part of the word problem, which is
always recursive.

 Now let us focus on the part of this procedure where Eve works with
 a particular presentation $\Gamma_i$ from Alice's pool. Suppose $\Gamma_i$ is not isomorphic to
 ${\hat \Gamma}_+$.   Since the
 ``no" part of the isomorphism problem between ${\hat \Gamma}_+$ and
$\Gamma_i$ is not recursive, Eve would have to   try out various
pairs of mappings between ${\hat \Gamma}_+$ and $\Gamma_i$ (see
above) indefinitely. Therefore, she will have to allocate
(indefinitely)  some memory resources to checking this particular
$\Gamma_i$. Since the number of $\Gamma_i$ grows exponentially with
the size of the presentation (which is the total length of
relators), Eve would require essentially unlimited storage space
and, in fact, she will reach physical limits (e.g. the number of
electrons in the universe) on the storage space very quickly
because, say, the number of presentations on 6 generators with the
total length of relators bounded by 100 is already more than
$10^{100}$.

 We note that there seems to be no way to bypass these storage space
requirements, even  assuming that Eve enjoys arbitrarily high (but
fixed)  computational speed. She basically has two options: (1) to
run in parallel subroutines corresponding to   each  $\Gamma_i$,
until one of the subroutines finds $\Gamma_i$ isomorphic to ${\hat
\Gamma}_+$; ~(2) to enumerate all steps in all subroutines the same
way that one enumerates rational numbers (so that Step $i$ in the
Subroutine  $j$  corresponds to the rational number $\frac{i}{j}$).

 The first option obviously requires essentially unlimited storage space. With the
second option, it may look like just having arbitrarily high
computational speed would be sufficient for Eve. However, since
under this arrangement Eve would have to interrupt each   subroutine
and then return to it, she would have to at least store the
information indicating where she left off  each  particular
subroutine. Therefore, since the number of subroutines is huge, she
would again reach physical limits on the storage space very quickly.

\section{Quotient attack}
\label{quotient}

In this section, we discuss an attack which is, in general, more
efficient (especially in real life) than the ``brute force" attack
described in Section \ref{isomorphism}. We use here some
group-theoretic terminology not supported by formal definitions when
we feel it should not affect the reader's understanding of the
material. Some of the basic terminology has to be introduced though.

A group $G$ is called {\it abelian} (or commutative)  if $[a, b] =1$
for any $a, b \in G$, where $[a, b]$ is the notation for $a^{-1}
b^{-1} ab$. Thus, $[a, b] =1$ is equivalent to $ab = ba$. This can
be generalized in different ways. A group $G$ is called {\it
metabelian} if $[[x, y], [z,t]] =1$ for any $x, y, z, t \in G$. A
group $G$ is called {\it nilpotent of class $c \ge 1$} if $[y_1,
y_2, \dots, y_{c+1}]=1$ for any $y_1, y_2, \dots, y_{c+1} \in G$,
where $[y_1, y_2, y_3]= [[y_1, y_2], y_3]$, etc.

We note that in the definition of an abelian group, it is sufficient
to require that $[x_i, x_j]=1$ for all {\it generators} $x_i, x_j$
of the group $G$. Thus, {\it any finitely generated abelian group
is finitely presented.} The same is true for all finitely generated
nilpotent groups of any class  $c \ge 1$, but not for all metabelian
groups. In particular, it is known that finitely generated {\it free
metabelian groups} are infinitely presented \cite{BST}. (A free
metabelian group is the factor group $F/[[F,F],[F,F]]$ of a free
group by the second commutator subgroup.)

Now we get to quotient attacks.
 One way for Eve to try to positively identify those places in Bob's binary sequence where
he intended to transmit a 0 is to use a {\it quotient test} (see
e.g. \cite{KMSS} for a general background). That means the
following: Eve tries to add finitely or infinitely many relators to
the given presentation ${\hat \Gamma}$ to obtain a presentation
defining a group $H$ with solvable word problem (more accurately, a
group $H$ that Eve can {\it recognize} as having solvable word
problem).

  It makes sense for Eve  to only try recognizable quotients, such as, for example,
 abelian or, more generally, nilpotent ones. This amounts to adding
 specific relators to ${\hat \Gamma}$; for example, for an abelian quotient, Eve
can add relators $[x_i', x_j']$ for all pairs of generators $x_i',
x_j'$ in ${\hat \Gamma}$. For nilpotent quotients, Eve will have to
add commutators of higher weight in the generators. For a metabelian
quotient, Eve will have to add infinitely many relators (because, as
we have already mentioned, free metabelian groups are infinitely
presented), but this is not a problem since she does not have to
``actually add" those relators; she can just consider ${\hat
\Gamma}$ as a presentation {\it in the variety of metabelian groups}
and apply the relevant algorithm for solving the word problem which
is universal for all groups finitely presented in the variety of
metabelian groups.

 Note that this trick will  {\it not} work with hyperbolic
 quotients, say. This is because there is no way, in general, to add
specific relators (finitely or infinitely many) to ${\hat \Gamma}$
to make sure that the extended presentation defines a hyperbolic
group. This deprives Eve from using a (rather powerful, cf.
\cite{KMSS}) hyperbolic quotient attack.

 Classes of groups with solvable word problem are summarized in the
 survey \cite{Olga}. It appears that a quotient attack can
 essentially employ either a nilpotent or a metabelian quotient of ${\hat
 \Gamma}$. This is why, to foil such attacks, Alice adds a relator  $x_i' = \prod_{j=1}^{M} [x_i', w_j]$
 to ${\hat \Gamma}$ (see our  Section \ref{pool}). This is also the
 reason why Bob should choose a   word of the form $[x_i', u]$ when
 he wants to transmit an element not equal to 1 in
$\Gamma'$ (see Section \ref{equal?}). Indeed, a metabelian quotient
attack on an element of the form $[x_i', u]$ will not
 work because this element belongs to the second commutator subgroup
 of the group defined by ${\hat \Gamma}$ since in this group, $x_i' = \prod_{j=1}^{M} [x_i', w_j]$,
so $x_i'$ belongs to the commutator subgroup of the given
 group.
  Furthermore, an element of the form $[x_i', u]$ belongs to
 every term of the lower central series of the given
 group since in this group, $[x_i', u]=[\prod_{j=1}^{M} [x_i', w_j], u]=
 [\prod_{j=1}^{M} [\prod_{j=1}^{M} [x_i', w_j], w_j], u]$, etc.
 This foils nilpotent quotient attacks, too.

\end{document}